# A Stochastic-Robust Approach to Hierarchical Generation-Transmission Expansion Planning

Maria de Luján Latorre, Gerson C. Oliveira, Ricardo C. Perez, Lucas Okamura and Silvio Binato
{lujan, gerson, ricardo, okamura, silvio}@psr-inc.com.
PSR
Rio de Janeiro, Brazil

*Abstract* – In this paper, the differences between an *integrated* and *hierarchical* generation and transmission expansion planning approaches are first described. Then, the hierarchical approach is described in detail. In general terms, in this scheme the investment decisions are made in steps, instead of an overall optimization scheme. This paper proposes a stochastic hierarchical generation-transmission expansion planning methodology based on a three-step procedure, as follows:
1. In this step, an integrated expansion planning problem of generation and interregional interconnections is solved;
2. Taking the optimal expansion plan of (generation and interconnections) into account, a production costing simulation with the detailed network representation is performed without monitoring circuit flow limits (except the ones in the interconnections, which are monitored); this simulation produces a set of optimal dispatch scenarios (vectors of bus loads and generation);
3. A Transmission Expansion Planning (TEP) model is then applied to determine the least-cost transmission expansion plan that is *robust* with respect to all operation scenarios of step 2, using an enhanced Benders decomposition scheme that: (a) incorporates a subset of the operation scenarios in the investment module; and (b) presents a "warm-up" step with a "greedy" algorithm that produces a (good) feasible solution and an initial set of feasibility cuts.

The application of the hierarchical planning scheme is illustrated with a realistic multi-country generation and transmission planning case study of the Central America's electricity market.

## I. NOMENCLATURE

| | |
|---|---|
| $s = 1, \ldots, S$ | operation scenarios |
| $j = 1, \ldots, J$ | existing circuits |
| $k = 1, \ldots, K$ | candidate circuits |
| $i = 1, \ldots, I$ | buses |
| $g^s$ | bus generation vector for scenario $s$ |
| $d^s$ | bus load vector for scenario $s$ |
| $\gamma$ | circuit susceptance vector |
| $c$ | vector of annualized investment cost (candidate circuits) |
| $\mu$ | (high) penalty cost for load curtailment |
| $A$ | incidence matrix (dimension $I \times (K+J)$), where the $n^{th}$ column has zeros in all rows, except the bus terminals of the $n^{th}$ circuit, $i_n$ and $j_n$ |
| $e$ | unitary vector $[1 \ \ldots \ 1]$ |
| $\bar{f}$ | vector of circuit flow limits |
| $r^s$ | bus load curtailment vector for scenario $s$ |
| $f^s$ | circuit flow vector for scenario $s$ |
| $\theta^s$ | bus voltage angle vector for scenario $s$ |
| $x$ | binary vector of circuit investment decisions |
| $M_k$ | "big M" constant |
| $p = 1, \ldots, P$ | indexes the Benders iterations |
| $\{\hat{x}_k^p\}$ | trial expansion plan produced by the investment module at the $p^{th}$ iteration |
| $S^p$ | subset of the operation scenarios that are (still) infeasible at the $p^{th}$ iteration |
| $\hat{w}_s^p$ | load curtailment resulting from trial expansion plan $\{\hat{x}_k^p\}$ and operating scenario $s$ |
| $\pi_{ks}^p$ | scalar variable that represents the greatest load curtailment |
| $\delta$ | Benders cut coefficient of investment variable $k$ for scenario $s$ at the $p^{th}$ iteration |

## II. INTRODUCTION

The origin of the expansion planning problem of electrical systems resides on the need for new investments in generation and transmission systems required to face the demand growth and meet planning criteria. In general terms, this decision process involves meeting economic, reliability and environmental criteria, within the framework of national policies on energy (such as decarbonization ones based on renewable penetration targets). In addition, one of the greatest challenges is how to deal with the uncertainties inherent in the planning process, such as the load growth, the hydrological inflows and the generation availability, especially in renewable-based systems. In that context, selecting the "best" of a group of alternatives is what characterizes the combinatorial nature of this expansion planning problem.

Taking all the aforementioned facts into account, this expansion planning problem can be modeled as a mixed-integer multistage stochastic problem. In this case, this will be an *integrated* generation and transmission (G&T) expansion planning approach. In [1], a methodology based on Benders decomposition and multistage stochastic optimization is proposed to solve this problem. The diagram below illustrates the main features of this methodology:

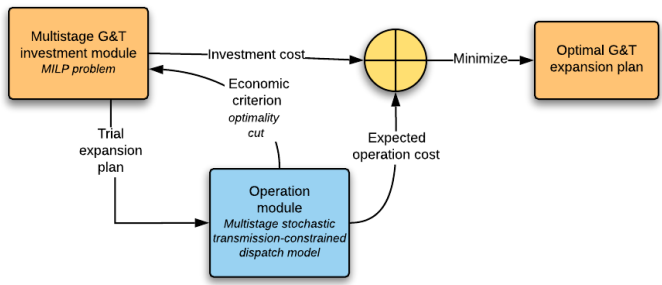

Fig. 1. Decomposition scheme for stochastic integrated generation and transmission expansion planning. Source: [1].

This methodology has been successfully applied to planning studies of several countries in Latin America [2], [3]. However, there are real-life situations in which the *integrated* planning approach may not be an option and a *hierarchical* approach, in which investment decisions are made in *steps*, instead of an overall optimization scheme is then applied [4], [5].

For example, since Brazil and Mexico are large-scale electric system, the number of G&T candidates and the consequent combinatorial nature of the G&T problem are huge. Consequently, the overall G&T optimization scheme is infeasible to be applied for real planning purposes.

Additionally, the hierarchical approach is also applied in some places due to country regulations. This usually occurs when there is a competitive environment in the generation sector (for example through auctions) and the transmission facilities must be decided (or anticipated from a long-term indicative generation-transmission expansion plan) based on the real information of the new generators' connections.

Finally, as will be shown among this document, the hierarchical approach is also adequate is planning for multi-country pools. In this paper, we describe a stochastic-robust programming methodology for generation and transmission planning of a multi-country system, and its application to the planning of Central America's electricity market.

The hierarchical approach is based on a three-step procedure, as follows:

- Step 1: the country is divided in regions and the network representation is simplified (usually by converting each electrical area into one single node and modeling only the maximum power transfer between regions). In the case of a multi-country system, each country is modeled as a single node having interconnections with the other countries. In this step, an integrated expansion planning problem of generation and interregional interconnections is solved;
- Step 2: taking the optimal expansion plan of (generation and interconnections) into account, a production costing simulation with the detailed network representation is performed without monitoring circuit flow limits (except the ones in the interconnections, which are monitored); this simulation produces a set of optimal dispatch scenarios (vectors of bus loads and generation);
- Step 3: a Transmission Expansion Planning (TEP)

model is then applied to determine the robust least-cost expansion plan that does not lead to overloads in any dispatch scenario; the optimization technique is based on Benders decomposition with feasibility cuts. An initial step with a heuristic "greedy" algorithm is used to determine an initial feasible solution and a set of strong Benders cuts.

The key issue in the TEP is how to deal with the uncertainties, which exist mainly due to the different hydrology in different river basins and strong volatility of the Variable Renewable Energy resources (VREs). There are TEP models that consider them based on Stochastic Programming (SP) [6], [7] and Robust Optimization [8], [9]. As can be seen, we use a SP approach in the generation and area interconnection planning, and a RO approach for the transmission planning of each area. For a more detailed TEP literature review, the reader should refer to [10], [11].

The next section details the hierarchical approach proposed in this paper.

### III. HIERARCHICAL GENERATION AND TRANSMISSION EXPANSION PLANNING

1. *Planning generation and interregional interconnections*

The first step is to determine the least-cost expansion plan of generation and interregional interconnections. This is done with the Benders decomposition/SDDP technique described in [1]. In this step a simplified network representation is applied modeling only the major transmission links between electrical areas.

The objective function of this optimization problem is to minimize investment (capital) and operation (fuel, O&M, deficit, etc.) costs. As can be seen, while modeling different regions (or countries) with interconnections (existing, planned and candidates), the main objective is evaluate the classic trade-offs "cheap and far" versus "more expensive, but near", i.e., make hydro plants and VREs, which are usually located far away from load centers compete against thermal plants which in turn are usually located near load centers. This concept may be also extended to the competition between the countries' energy mix and available resources in the case of a multi-country system.

2. *Preparation of load/generation scenarios and proposition of transmission candidates for the TEP*

Once the optimal generation-interconnection expansion plan is obtained, the network representation is replaced to the (real) detailed one and a production costing simulation is performed respecting circuit flow limits only in the interregional interconnections. This probabilistic simulation provides a set of bus generation and load vectors, known as operation scenarios, that will be used in the TEP described next. Besides that, the results of this dispatch simulation are used by the transmission planner in order to deeply evaluate the system loading and propose transmission candidates (including new rights-of-way when feasible).

3. *Robust transmission expansion planning*

The objective of this step in the hierarchical planning process is to determine a least-cost transmission network reinforcement that is robust with respect to all dispatch scenarios, i.e., that has no overloads in any scenario.

This problem is solved by a second Benders decomposition scheme summarized below:

- The investment module determines the transmission reinforcements;
- The operation modules are linearized optimal power flow models (DC OPFs) that minimize the load shedding required to eliminate all overloads; and
- The feedback from the DC OPFs to the investment module are Benders feasibility cuts.

In section IV, we describe in more detail the optimization methodology for transmission planning under uncertainty with emphasis on two improvements which are the main contributions of this paper and led to substantial reductions in computational effort:

- Representation of a set of operation scenarios in the investment module; and
- A two-step solution scheme in which a "greedy" algorithm is used first to obtain a good feasible solution plus a set of feasibility cuts that serve as a "hot-start" for the Benders-based second step.

In section V, we describe the application of the hierarchical planning scheme to a realistic study for Central America's MER.

## IV. ROBUST TRANSMISSION EXPANSION PLANNING WITH MULTIPLE OPERATION SCENARIOS

In this section, we describe the optimization methodology for transmission planning under uncertainty. For ease of presentation, we initially formulate the transmission planning problem as a single (and huge) Mixed-Integer Linear Programming (MILP) problem for a given year, which in turn is a DC Optimum Power Flow (DC OPF) model. Next, we describe the Benders decomposition scheme, followed by enhanced investment module and, finally, the first-step "greedy" algorithm.

### A. Transmission Planning Formulation as one MILP Problem

#### 1. Basic problem formulation

*Objective function*

The objective is to minimize the sum of investment costs plus a penalty for infeasibilities (load curtailments due to overloads) in the operation scenarios.

$$Min \sum_k c_k x_k + \mu \times \sum_s e' r^s \qquad (1a)$$

*Power balance in each bus (Kirchhoff's first law)*

$$Af^s + r^s = d^s - g^s \qquad \forall s \qquad (1b)$$

*Flow equation for existing circuits (Kirchhoff's second law)*

$$f_n^s = \gamma_n \Delta\theta_n^s \quad (\Delta\theta_n^s = [\theta_{i_n}^s - \theta_{j_n}^s]) \qquad \forall n, s \qquad (1c)$$

*Flow equation for candidate circuits (disjunctive formulation)*

$$0 \leq |f_k^s - \gamma_k \Delta\theta_k^s| \leq M_k(1 - x_k) \qquad \forall k, s \qquad (1d)$$

As usual when disjunctive constraints are used, and discussed in [12], for the case of transmission planning, the value of this constant is critical for the efficient solution of the MILP problem. [13], [12] and show how to calculate the smallest possible "big M" for each right-of-way.

*Circuit flow limits*

$$|f_j^s| \leq \bar{f}_j \qquad \forall j, s \qquad (1e)$$

$$|f_k^s| \leq \bar{f}_k x_k \qquad \forall k, s \qquad (1f)$$

*Limits on load curtailment*

$$r^s \leq d^s \qquad \forall s \qquad (1g)$$

*Binary investment variables*

$$x_k \in \{0,1\} \qquad \forall k \qquad (1h)$$

It can be seen that although the problem size increases linearly with the number of dispatch scenarios, the number of binary variables $x_k$ remains the same. Additionally, it is worth noting that these binary variables that are responsible for "coupling" OPF equations of the different dispatch scenarios, i.e., the optimal expansion plan must meet all of them.

Furthermore, the Formulation (1) is strengthened by three sets of additional constraints, described next:

- Export/import bottlenecks;
- Disconnected paths; and
- Precedence constraints for duplicate circuits.

#### 2. Export/Import bottlenecks

Let $\delta_m^s = |g_m^s - d_m^s|$ be the absolute value of the total power entering/leaving bus $m$ in operating scenario $s$. Let $\delta_m^* = Max_s\{\delta_m^s\}$. Denote by $K_m$ the set of candidate circuits directly connected to bus $m$, and let $J_m$ denote the set of exiting circuits directly connected to bus $m$. If the total capacity of the existing circuits directly connected to $m$ is smaller than $\delta_m^*$, there is a "power flow bottleneck" in the bus and the following feasibility constraint is added:

$$\sum_{k \in K_m} \bar{f}_k x_k \geq \delta_m^* - \sum_{n \in J_m} \bar{f}_n \qquad (2)$$

#### 3. Disconnected paths

A new expansion path may be composed of several circuits in series, which start and end in buses connected to the main grid, but that have disconnected intermediate "passage" buses (new candidate substations), which are also called transshipment buses. Let $n$ be an intermediate bus, and $D_n$ the set of candidate circuits directly connected to it. The following

constraints ensure that the intermediate bus either remains disconnected or becomes part of a connected path:

$$x_i \leq \sum_{j \in D_n, j \neq i} x_j \qquad i \in D_n \qquad (3)$$

For example, suppose that only two candidate circuits, $(i,n)$ and $(n,j)$, are connected to candidate bus $n$. The two corresponding constraints (3) imply that $x_{(i,n)} = x_{(n,j)}$.

Although the connection constraints (3) are applied to each candidate bus, they automatically take into account the situation of two candidate buses $i$ and $j$ that are the terminals of a candidate circuit $(i,j)$, i.e. ensure that either both buses remain disconnected or that a candidate path will be built, composed of this circuit and at least two more candidate circuits starting and ending in connected buses. The same logic applies to candidate paths formed by more than two candidate circuits; and to two or more candidate paths having in common a candidate bus $n$ and a candidate circuit $(i,n)$, e.g., paths $i \to j \to n$ and $i \to n \to k$.

*4. Precedence constraints*

In order to avoid symmetry due to $w$ multiple equal candidates connected to terminal buses $(i,j)$ the following logical constraint is written:

$$x_{(i,j)}^{k+1} \leq x_{(i,j)}^{k} \qquad \forall k < w \qquad (4)$$

Avoiding the so called symmetry effect is important because it slows down the performance of the Branch-and-Bound algorithm, caused by variables with the same costs in the objective function and the same coefficients in the constraints. Furthermore, as could be seen, the extension of the above precedence constraint to more than two candidates is trivial, each inequality is written for a pair of consecutively numbered candidates in the same corridor $k$.

### B. Decomposition Scheme

Due to the typically large number of operation scenarios, the direct MILP formulation (1) is computationally infeasible for realistic planning studies. The Benders decomposition scheme mentioned previously is attractive because it allows the application of an iterative solution approach: one sub-problem (master), the *investment module*, is related to the investment decisions, and solved by a MILP problem; and the other sub-problems (slaves), the *operation modules*, are related to the feasibility evaluation of each dispatch scenario, and solved by LP techniques. The decomposition scheme is summarized in the figure below:

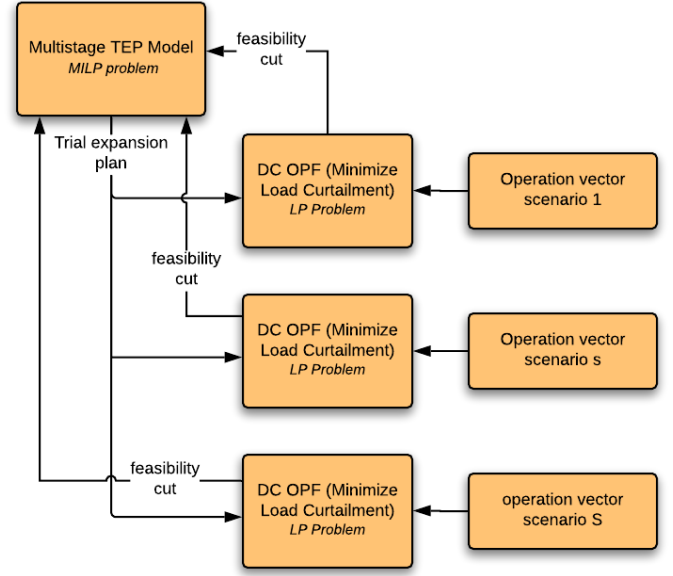

Fig. 2. Decomposition-based transmission planning under uncertainty.

*1. Basic investment module formulation*

The "orthodox" formulation of the investment module would be:

$$Min \sum_k c_k x_k + \mu \times \delta \qquad (5a)$$

$$\widehat{w}_s^p + \sum_k \pi_{ks}^p (x_k - \hat{x}_k^p) - \delta \leq 0 \qquad \forall p, s \in S^p \qquad (5b)$$

*2. Feasibility sub-problem formulation*

Given a trial investment plan $\{\hat{x}_k^p\}$ and an operating scenario $s \in S^p$ in the $p^{th}$ iteration of the Benders decomposition scheme, the feasibility assessment sub-problem is formulated as:

$$w_s^p = Min \sum_i r_i^s \qquad (6a)$$

$$Af^s + r^s = d^s - g^s \qquad \forall s \qquad (6b)$$

$$f_j^s = \gamma_j \Delta \theta_j^s \qquad \forall j \qquad (6c)$$

$$0 \leq |f_k^s - \gamma_k \Delta \theta_k^s| \leq M_k^p \; ( \\ = M_k(1 - \hat{x}_k^p)) \qquad \forall k \; \leftarrow \pi_k^M \qquad (6d)$$

$$|f_j^s| \leq \bar{f}_j \qquad \forall j \qquad (6e)$$

$$|f_k^s| \leq \bar{f}_k^p \; (= \bar{f}_k \hat{x}_k^p) \qquad \forall k \; \leftarrow \pi_k^f \qquad (6f)$$

$$r^s \leq d^s \qquad (6g)$$

3. *Benders cut calculation*

The Benders cut coefficient $\pi_{ks}^p$ in equation (5b) is calculated as:

$$\pi_{ks}^p = -M_k \pi_k^M + \bar{f}_k \pi_k^f \qquad \forall k \qquad (7)$$

Where $\pi_k^M$ and $\pi_k^f$ are the simplex multipliers associated to constraints (6d) and (6f) at the optimal solution of sub-problem (6).

## C. Enhanced Decomposition Scheme

As mentioned, the Benders formulation (4)-(6) may converge slowly for some larger systems. For this reason, we developed two improvements to the procedure: (i) incorporate a set of operation scenarios to the investment module; and (ii) first apply a greedy algorithm to obtain a feasible solution plus feasibility cuts to "hot-start" the Benders scheme, which lead to formulation presented below.

1. *Enhanced investment module*

The enhanced investment module is formulated as follows:

$$Min \sum_k c_k x_k + \mu \times \sum_{s \in S^c} \sum_i r_i^s + \mu \times \delta \qquad (8a)$$

$$Af^s + r^s = d^s - g^s \qquad \forall s \in S^c \qquad (8b)$$

$$f_j^s = \gamma_j \Delta \theta_j^s \qquad \forall j, s \in S^c \qquad (8c)$$

$$0 \leq |f_k^s - \gamma_k \Delta \theta_k^s| \leq M_k(1 - x_k) \qquad \forall k, s \in S^c \qquad (8d)$$

$$|f_j^s| \leq \bar{f}_j \qquad \forall j, s \in S^c \qquad (8e)$$

$$|f_k^s| \leq \bar{f}_k x_k \qquad \forall k, s \in S^c \qquad (8f)$$

$$r^s \leq d^s \qquad \forall s \in S^c \qquad (8g)$$

$$\hat{w}_s^p + \sum_k \pi_{ks}^p (x_k - \hat{x}_k^p) - \delta \leq 0 \qquad \forall p, s \in S^p \qquad (8h)$$

The subset $S^c$ is composed of the most severe scenarios (also called *critical* scenarios), selected as follows.

2. *Selection of the critical scenarios for the investment module*

a) Define the critical scenario selection criterion: greater value of load shedding or greater distribution of load shedding (two different measures of the *severity* of the transmission problems for the scenario $s$).
b) Define the subset $S^c$ (user-defined parameter) which is the maximum number of scenarios that can be contemplated in the master problem.
c) Run a linearized power flow model for all dispatch scenarios $s = 1, ..., S$ and select the subset $S^0$ of scenarios that have at least one circuit overload, also called *severe* scenarios.
d) For each scenario $s \in S^0$, run a feasibility problem (5) and calculate $w_s^0$ and $i_s^0$, which are the minimum amount of load shedding in the entire system required to eliminate the overloads and the number of buses with any load shedding, respectively.
e) Depending on the options selected at a) and b), rank the scenarios $s \in S^0$ by decreasing $w_s^0$ or $i_s^0$ and select the first $S^c$ to be added to the investment problem.

## D. Initial "Greedy" Algorithm

The "greedy" algorithm uses the same idea of adding a subset of operating scenarios to the MILP investment problem, but applies the procedure sequentially:

(i) Apply the procedure (a)-(c) to find the subsets $S^0$ and $S^c$ of scenarios.
(ii) Use the MILP formulation (1)-(4) to determine the optimal plan for the $S^c$ scenarios selected in step (i) and add the resulting reinforcements to the network.
(iii) Given the reinforced network, verify if there are remaining overloads in the remaining operation scenarios ($S^0 - S^c$). This is done by running the linearized power flow model for each scenario. If there are overloads, go to step (iv). Otherwise, go to step (v).
(iv) Select the new $S^c$ critical scenarios of step (iv) and return to step (ii) to solve the MILP formulation (with the reinforced network).
(v) In this final step, we model the feasibility sub-problems in order to remove circuits that became redundant due to the sequential (i.e. "myopic") nature of the greedy scheme. This is done by ranking the added circuits by decreasing investment cost and eliminating them if their removal does not cause any overload in all operation scenarios. In contrary, if a circuit removal causes any overload, then, a feasibility cut (7) is calculated. As mentioned, those feasibility cuts will be used later as a "hot-start" for the enhanced Benders decomposition scheme (8).

## E. Final Two-step Transmission Planning Algorithm

As mentioned, we start with the greedy scheme (i)-(vi), followed by the enhanced decomposition scheme (8), "warmed" by an initial set of feasibility cuts from step (iv) of the greedy procedure. Furthermore, another byproduct of the feasible solution is an upper bound to the Benders decomposition. However, our computational experience showed no improvement with this measure.

It is worth noting that the algorithm proposed in this paper is fully applicable to contingency analysis and N-1 planning. In general, the algorithm described in this section is first applied to find the base case optimal expansion plan and then the whole process starts again for contingency assessment and N-1 expansion planning.

Finally, before moving the case study, the figure presented below summarizes the stochastic hierarchical planning approach proposed in this paper.

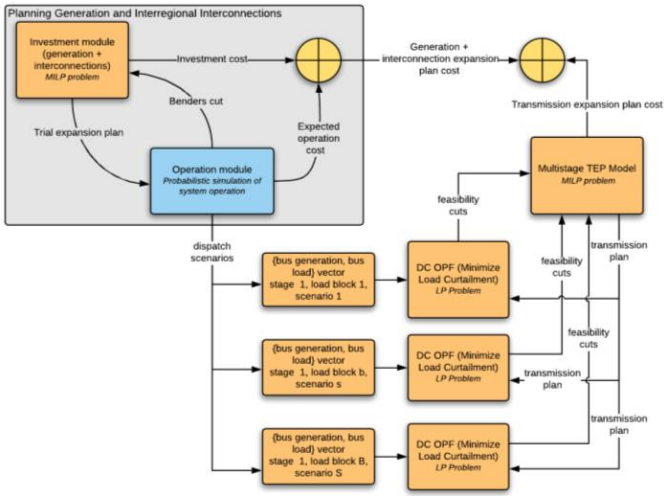

Fig. 3. The stochastic hierarchical planning approach proposed in this paper.

## V. CASE STUDY

### A. The Regional Electricity Market (MER)

The MER is currently composed of six countries: Panama, Costa Rica, Nicaragua, Honduras, El Salvador and Guatemala.

The figure below shows the main characteristics of each country (installed capacity). We see that there is a wide mix of generation technologies, with a historically strong hydro share and, more recently, a fast penetration of wind, solar and biomass.

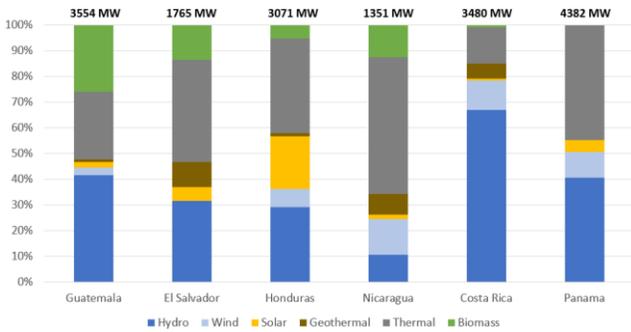

Fig. 4. Installed capacity of the MER [14].

Additionally, the figure below shows the current MER interconnections and their respective capacities. We also show the interconnection between Mexico and Guatemala, which is not part of the planning study because Mexico does not belong to MER (the power exchanges between Mexico and Guatemala are represented as operation scenarios).

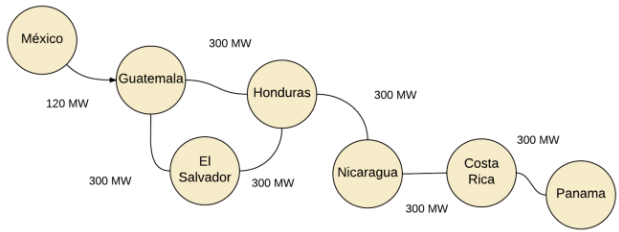

Fig. 5. MER interconnection (plus Mexico-Guatemala)

The following figure shows the MER existing transmission network represented in the planning study.

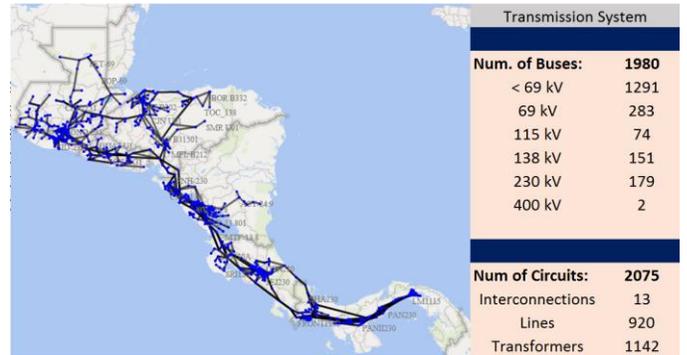

Fig. 6. MER transmission network general overview.

### B. Probabilistic Simulation of System Operation before Reinforcements

The study planning horizon of this case study is: 2020-2024. In the previous section, we showed the current installed capacity and generation mix of each MER country. The figure below presents the annual demand forecasts (GWh) for the study horizon.

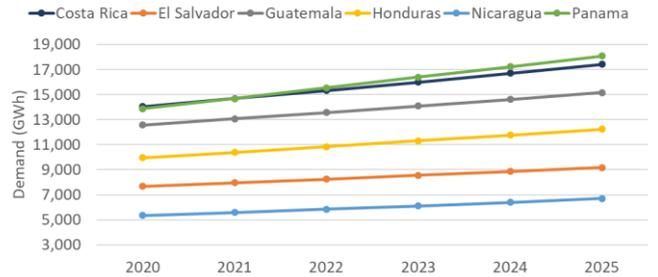

Fig. 7. Demand forecasts for the MER planning study

Having the demand forecasts, the first step was to carry out a probabilistic simulation of system operation representing both generation and the transmission network, including only the reinforcements already under construction (the SDDP model [15] was used for this task). This simulation of has monthly time stages, with 5 load blocks (of the Load Duration Curve representation) in each month and 25 renewable energy scenarios (inflows for the hydro plants and generation scenarios for VREs maintaining the temporal and spatial correlations).

Table I shows the expected value of system short-run marginal costs (SRMCs) and the deficit risks for each year. The systemic SRMC of each *configuration* (where configuration is a combination of month, renewable energy scenario and load block) is the weighted average of the SRMCs over all buses (the weights are the bus loads). In turn, the expected annual SRMC is the average of the systemic SRMCs over all hydrology/renewable scenarios, monthly stages and load blocks in each stage. The deficit risk indicates which percentage of scenarios present any deficit of any magnitude in at least one configuration.

TABLE I
MER: EXPECTED ANNUAL SRMCs AND DEFICIT RISKS (NO REINFORCEMENTS)

| Year | Panama ($/MWh) | (%) | Costa Rica ($/MWh) | (%) | Nicaragua ($/MWh) | (%) | Honduras ($/MWh) | (%) | El Salvador ($/MWh) | (%) | Guatemala ($/MWh) | (%) |
|---|---|---|---|---|---|---|---|---|---|---|---|---|
| 2020 | 134 | 2 | 133 | 2 | 138 | 2 | 133 | 2 | 128 | 0 | 127 | 0 |
| 2021 | 154 | 11 | 154 | 5 | 156 | 11 | 151 | 26 | 142 | 0 | 141 | 0 |
| 2022 | 188 | 93 | 186 | 86 | 187 | 100 | 180 | 100 | 157 | 0 | 156 | 0 |
| 2023 | 265 | 100 | 261 | 100 | 254 | 100 | 241 | 100 | 183 | 3 | 183 | 6 |
| 2024 | 386 | 100 | 376 | 100 | 358 | 100 | 341 | 100 | 205 | 9 | 203 | 11 |

As expected, in later years the expected annual SMRCs increase to very high values as the deficit risks, signaling the need for reinforcements.

### C. Planning of Generation and Regional Interconnections

The next step is to determine the optimal expansion of generation and regional interconnections. The methodology is outlined in Figure 1 and described in detail in reference [1]. Basically, it is a Benders decomposition scheme where the investment module decides the capacity expansion of generation and regional interconnections, and the stochastic operation model is solved by the SDDP model as illustrated in Figure 1. Tables II and III present the number of candidate generators and interconnections considered in this step, respectively.

TABLE II
CANDIDATE GENERATION PROJECTS IN EACH COUNTRY

| Plant Type | Costa Rica | | El Salvador | | Guatemala | |
|---|---|---|---|---|---|---|
| | Num. | Total Cap. (MW) | Num. | Total Cap. (MW) | Num. | Total Cap. (MW) |
| Hydro | 3 | 709 | 3 | 397 | 1 | 60 |
| Diesel | 6 | 700 | 6 | 700 | 6 | 700 |
| Coal | 6 | 1,800 | 6 | 1,800 | 6 | 1,800 |
| Bunker | 2 | 200 | 2 | 200 | 2 | 200 |
| Natural Gas | 2 | 1,000 | 3 | 1,000 | 2 | 1,000 |
| Geothermal | 2 | 110 | 1 | 30 | 3 | 77 |
| Wind | 3 | 300 | 3 | 300 | - | - |
| Solar | - | - | - | - | - | - |
| Biomass | - | - | - | - | - | - |

| Plant Type | Honduras | | Nicaragua | | Panama | |
|---|---|---|---|---|---|---|
| | Num. | Total Cap. (MW) | Num. | Total Cap. (MW) | Num. | Total Cap. (MW) |
| Hydro | 7 | 911 | 10 | 767 | 11 | 354 |
| Diesel | 6 | 700 | 6 | 700 | 6 | 700 |
| Coal | 6 | 1,800 | 6 | 1,800 | 6 | 1,800 |
| Bunker | 2 | 200 | 2 | 200 | 2 | 200 |
| Natural Gas | 3 | 1,000 | 3 | 1,000 | 5 | 1,310 |
| Geothermal | - | - | 8 | 230 | - | - |
| Wind | 5 | 395 | 4 | 320 | 7 | 416 |
| Solar | - | - | - | - | - | - |
| Biomass | - | - | 1 | 30 | - | - |

TABLE III
CANDIDATE INTERCONNECTION PROJECTS

| Interconnection | Num. | Total Cap. (MW) |
|---|---|---|
| Nicaragua-Costa Rica | 5 | 1,732 |
| Honduras - Nicaragua | 5 | 1,699 |
| Panama - Costa Rica | 5 | 1,233 |
| Honduras - El Salvador | 3 | 1,008 |
| Guatemala - Honduras | 1 | 374 |
| Guatemala - El Salvador | 1 | 374 |

Th task has been run on a virtual server on AWS with 16 CPUs. The numbers of Benders iterations (solving the investment problem MILP plus a SDDP-based stochastic operation model) was 14; Since these MILP problems are small, we use a very tight MILP convergence tolerance of 0.001%; the total wall clock time was two hours. Tables IV and V show the optimal expansion plan of this step.

TABLE IV
OPTIMAL GENERATION EXPANSION PLAN OF THE MULTI-COUNTRY SYSTEM

| Year | Coal (MW) | Hydro (MW) | Wind (MW) | Geothermal (MW) | Biomass (MW) |
|---|---|---|---|---|---|
| 2020 | 0 | 290 | 114 | 0 | 30 |
| 2021 | 0 | 253 | 95 | 20 | 0 |
| 2022 | 500 | 147 | 0 | 0 | 0 |
| 2023 | 0 | 261 | 0 | 125 | 0 |
| 2024 | 0 | 420 | 0 | 85 | 0 |

TABLE V
OPTIMAL INTERCONNECTION EXPANSION PLAN OF THE MULTI-COUNTRY SYSTEM

| Bus From Name | Voltage (kV) | System From | Bus To Name | Voltage (kV) | System To | Circuit Cost (k$) | Circuit Rating (MW) | Entry Year |
|---|---|---|---|---|---|---|---|---|
| PRD B618 | 230 | Honduras | FNH-230 | 230 | Nicaragua | 4,698 | 317 | 2020 |
| FRONTPRO | 230 | Panama | RCL230A | 230 | Costa Rica | 5,617 | 317 | 2020 |
| FRONTCHA | 230 | Panama | CAH230 | 230 | Costa Rica | 7,674 | 200 | 2020 |
| FRONTPRO | 230 | Panama | RCL230A | 230 | Costa Rica | 5,617 | 259 | 2020 |
| FRONTPRO | 230 | Panama | RCL230A | 230 | Costa Rica | 5,617 | 200 | 2020 |
| FRONTDOM | 230 | Panama | RCL230B | 230 | Costa Rica | 4,186 | 200 | 2020 |
| NAC 230 | 230 | Honduras | F-15SE-AG1 | 230 | El Salvador | 4,849 | 374 | 2021 |

### D. Transmission Expansion Planning for each Country

As described in the previous section, the optimal solution of the previous step produces operation scenarios (bus generation/load vectors) for each *configuration* (month, hydro/renewable generation scenario and load block). As can be seen, for the transmission planning study, which is performed per year (forward in time), the total number of operation scenarios per year is: 12 (months) × 5 (load blocks) × 25 (inflows + VREs), i.e., 1500 vectors of bus loads and generation setpoints.

The TEP model described in section IV has been applied to determine the robust least-cost expansion plan (no overloads for all 1500 scenarios per year) in each country (separately in order to reduce the size of the MILP problems being solved). In this case study, only base case is under evaluation.

Because the planning studies were independent (contemplating the imports/exports as injections at the frontier buses), they were carried out in parallel in the same 16-CPU server used for the generation/regional interconnection plan of the previous section. The list of candidate circuits considered in this step was composed by 183 transmission lines and 11 transformers. Since the MILP problems being solved in this task are complex (mainly due to the disjunctive constraints), we use MILP convergence tolerance of 3% instead of the tighter previous one (0.001%).

Table VI presents the optimal transmission expansion plan of the entire region:

TABLE VI
OPTIMAL TRANSMISSION EXPANSION PLAN (ALL COUNTRIES)

| Bus From Name | Bus From Voltage (kV) | Bus To Name | Voltage (kV) | System | Circuit Cost (k$) | Circuit Rating (MW) | Entry Year |
|---|---|---|---|---|---|---|---|
| PAR230 | 230 | PAL230B | 230 | Costa Rica | 15,737 | 374 | 2020 |
| PRO230 | 230 | BOQIII230 | 230 | Panama | 5,205 | 556 | 2020 |
| FRONTPRO | 230 | PRO230 | 230 | Panama | 2,832 | 193 | 2020 |
| MOI230A | 230 | CAH230 | 230 | Costa Rica | 7,408 | 259 | 2020 |
| RCL230A | 230 | PAL230A | 230 | Costa Rica | 8,272 | 300 | 2020 |
| CHA230 | 230 | FRONTCHA | 230 | Panama | 3,559 | 304 | 2020 |
| CAR230A | 230 | GEN230 | 230 | Costa Rica | 12,297 | 659 | 2020 |
| PAR230 | 230 | PAL230B | 230 | Costa Rica | 19,868 | 374 | 2020 |
| FRONTPRO | 230 | PRO230 | 230 | Panama | 2,832 | 193 | 2020 |
| FRONTPRO | 230 | PRO230 | 230 | Panama | 2,832 | 193 | 2020 |
| RCL230A | 230 | PAL230A | 230 | Costa Rica | 8,272 | 300 | 2020 |
| CRL B501 | 138 | RLN B521 | 138 | Honduras | 4,500 | 152 | 2021 |
| RLN B521 | 138 | CAR B540 | 138 | Honduras | 4,500 | 152 | 2021 |
| VNU B520 | 138 | CAR B540 | 138 | Honduras | 4,500 | 152 | 2021 |
| SFR115 | 115 | CBA115 | 115 | Panama | 1,202 | 230 | 2022 |
| LBS-138 | 138 | NAG-138 | 138 | Nicaragua | 4,097 | 80 | 2023 |
| AHUA-115 | 115 | SONS-115 | 115 | El Salvador | 3,696 | 130 | 2024 |
| RGU B518 | 138 | SIS B548 | 138 | Honduras | 7,607 | 152 | 2024 |
| GUA138 | 138 | FIL138 | 138 | Costa Rica | 4,070 | 265 | 2024 |
| CTE B513 | 138 | SIS B548 | 138 | Honduras | 1,474 | 152 | 2024 |

Table VII shows the total CPU time (in seconds) for the transmission expansion plan determination of each country in each year applying all steps described in section IV. Results are shown only for the cases where the TEP model needed to be executed, i.e., *severe* scenarios were identified during the linearized power flow executions for all configurations (task which took an average of 89 seconds to perform per country per year).

TABLE VII
EXECUTION TIME FOR THE TEP OF EACH COUNTRY

| System | 2020 plan (s) | 2021 plan (s) | 2022 plan (s) | 2023 plan (s) | 2024 plan (s) |
|---|---|---|---|---|---|
| Costa Rica | 291.05 | | | | 236.78 |
| El Salvador | | | | | 137.33 |
| Honduras | | | 221.16 | | 218.94 |
| Nicaragua | | | | 162.55 | |
| Panama | 666.64 | | 248.38 | | |

Given that each country's expansion is carried out separately in parallel, the wall-clock time is the highest sum along all years.

*E. Probabilistic Simulation of System Operation after the Reinforcements*

Finally, an SDDP-based probabilistic simulation of system operation representing both generation and transmission was carried out with the planned regional and country-level reinforcements. Table VIII shows the expected system SRMCs and deficit risks for each country and each year.

TABLE VIII
SYSTEM SRMCs AFTER TRANSMISSION REINFORCEMENTS

| Year | Panama | | Costa Rica | | Nicaragua | | Honduras | | El Salvador | | Guatemala | |
|---|---|---|---|---|---|---|---|---|---|---|---|---|
| | ($/MWh) | (%) | ($/MWh) | (%) | ($/MWh) | (%) | ($/MWh) | (%) | ($/MWh) | (%) | ($/MWh) | (%) |
| 2020 | 99 | 0 | 99 | 0 | 99 | 0 | 98 | 0 | 96 | 0 | 97 | 0 |
| 2021 | 92 | 0 | 92 | 0 | 92 | 0 | 91 | 0 | 91 | 0 | 91 | 0 |
| 2022 | 88 | 0 | 88 | 0 | 89 | 0 | 89 | 0 | 89 | 0 | 89 | 0 |
| 2023 | 95 | 0 | 95 | 0 | 95 | 0 | 95 | 0 | 95 | 0 | 95 | 0 |
| 2024 | 94 | 0 | 94 | 0 | 94 | 0 | 93 | 0 | 93 | 0 | 93 | 0 |

As expected, besides presenting no deficit risk, the multi-country system faces lower annual average SRMCs with representatively smaller differences between countries, showing that the G&T expansion plan has been optimized enabling though the interconnections to take advantage of the regional portfolio and available resources in each country to minimize investments and operating costs for the entire region.

## VI. CONCLUSIONS

This paper proposes a new stochastic-robust approach to hierarchical generation-transmission expansion planning. The main contributions of this work are in the TEP task under uncertainty which led to substantial reductions in computational effort:

- Representation of a set of operation scenarios in the investment module; and
- A two-step solution scheme in which a "greedy" algorithm is used first to obtain a good feasible solution plus a set of feasibility cuts that serve as a "hot-start" for the Benders-based second step.

This paper shows that the proposed methodology is able to find the robust least-cost expansion plan contemplating a wide range of dispatch scenarios while expanding the electrical network, which is a key feature for hydrothermal systems (due to the hydrology associated to the rivers located in different regions) and for the worldwide increasing penetration of VREs (due to the uncertainty and strong volatility associated).